\documentclass[11pt]{amsart}

\usepackage[T1]{fontenc}
\usepackage{lmodern}
\usepackage{amsmath,amssymb,amsthm}
\usepackage{xcolor}
\usepackage{graphicx}
\usepackage{mathdots}
\usepackage{hyperref}

\definecolor{darkblue}{rgb}{0,0,0.7}
\hypersetup{colorlinks=true,citecolor=darkblue,linkcolor=darkblue,urlcolor=darkblue}

\newcommand\defn[1]{\emph{\textcolor{darkblue}{#1}}}

\newtheorem{thm}{Theorem}[section]
\newtheorem{lemma}[thm]{Lemma}
\theoremstyle{definition}
\newtheorem{ex}[thm]{Example}

\numberwithin{figure}{section}
\numberwithin{equation}{section}

\newcommand{\R}{\mathbb{R}}
\newcommand{\Z}{\mathbb{Z}}
\newcommand{\N}{\mathbb{N}}
\newcommand{\C}{\mathbb{C}}

\newcommand{\la}{\lambda}

\newcommand{\rQ}{\mathrm{Q}}
\newcommand{\Ehr}{{E}}
\newcommand{\GT}{\mathrm{GT}}

\renewcommand{\a}{\mathbf{a}}
\renewcommand{\b}{\bb}
\newcommand{\x}{\bx}
\newcommand{\f}{\textbf{f}}

\newcommand{\ts}{\hskip0.03cm}
\renewcommand{\.}{\hskip0.06cm}
\newcommand{\nts}{\hskip-0.03cm}

\def\zz{\mathbb Z}

\def\nn{\mathbb N}

\def\cc{\mathbb C}

\def\la{\lambda}
\def\ga{\gamma}

\def\<{\langle}
\def\>{\rangle}

\def\Z{\zz}

\def\GL{ {\text {\rm GL} } }

\def\0{{\mathbf 0}}

\def\.{\hskip.06cm}
\def\ts{\hskip.03cm}

\def\bx{{\textbf{x}}}

\def\nts{\hskip-.04cm}

\def\bb{\mathbf{b}}

\def\C{C}

\begin{document}

\title[On stretched Kostka polynomials]{On stretched Kostka polynomials}

\author[Alexey Garber, Igor Pak, and Zachary Slonim]{Alexey Garber$^\star$, Igor Pak$^\diamond$, and Zachary Slonim$^\diamond$}

\thanks{\emph{Date:} \today}

\thanks{\thinspace${\hspace{-.45ex}}^\star$School of Mathematical \& Statistical Sciences, University of Texas Rio Grande Valley, Brownsville, TX 78520.\. Email:  \texttt{alexey.garber@utrgv.edu}}

\thanks{\thinspace${\hspace{-.45ex}}^\diamond$Department of Mathematics, UCLA, Los Angeles, CA 90095.\. Email: \texttt{\{pak,zachslonim\}@math.ucla.edu}}

\begin{abstract}
Stretched Kostka polynomials are defined by \ts$f_{\la,\mu}(N) := K_{N\la,N\mu}$\ts, where $K_{\la,\mu}$ denotes the corresponding Kostka number. In \cite{KTT04}, King--Tollu--Toumazet conjectured that, for every integer $k$, there are only finitely many polynomials $f_{\la,\mu}(N)$ such that $f_{\la,\mu}(1)=k$. For $k=3$, KTT also conjectured that there are exactly two such polynomials.
We prove both conjectures via Ehrhart theory, by counting lattice points in families of rational polytopes with a fixed set of normals.
\end{abstract}

\maketitle

\section{Introduction}

\subsection{Main results}

\defn{Kostka numbers} \ts $K_{\la,\ts\mu}$ \ts are ubiquitous in algebraic combinatorics and representation theory.
Notably, they can be defined as \emph{weight multiplicities} in an irreducible representation of $\GL(n,\cc)$.
Equivalently, they can be defined as transition coefficients between Schur and monomial symmetric functions:
$$
s_\la \ = \ \sum_{\mu} \. K_{\la,\ts\mu} \. m_\mu\,.
$$
Kostka numbers have several combinatorial interpretations, as the number of semistandard Young tableaux (SSYT) and
as the number of integer points in \emph{Gelfand--Tsetlin polytopes} \cite{GT50,GZ85}.  We refer to \cite{Ful97,Mac95,Sta99}
for the introduction and some background.

For an integer $N\ge 1$, write $N\la:=(N\la_1,N\la_2,\ldots)$ and define the \defn{stretched Kostka polynomial}
\[f_{\la,\ts\mu}(N)\ :=\ K_{N\la,\ts N\mu} \,.
\]
Famously, the polynomiality follows from the fermionic formula of Kirillov and Reshetikhin \cite[Thm.~4.2]{KR88},
see also \cite[Cor.~4.4]{Kir01}. In fact, polynomiality in $N$ is a remarkable regularity property
known in much greater generality,
see e.g.\ \cite{BGR04,McA08}. Geometrically, $f_{\lambda,\ts\mu}$ counts lattice points in dilations
of a Gelfand--Tsetlin polytope. In this setting, it coincides with the \emph{Ehrhart polynomial},
see e.g. \cite{Bar08, BR15}.

Motivated by this connection, in 2004, King, Tollu, and Toumazet (KTT) stated several conjectures
on the rigidity of these polynomials.  Notably, they conjectured in \cite[Conj.~3.3]{KTT04},
that there are only finitely many stretched Kostka polynomials $f_{\la,\ts\mu}(N)$ with
$K_{\la,\ts\mu}=k$ for a fixed integer~$k$.  Our main result is a complete resolution
of this conjecture:

\smallskip

\begin{thm}[{\rm formerly KTT's conjecture}{}]\label{t:finite}
Fix $k\in\N$. Then the following family of polynomials is finite:
\[
P_k\ :=\ \bigl\{\ts f_{\la,\ts\mu}(N)\, : \, f_{\la,\ts\mu}(1)=k\ts\bigr\}.
\]
\end{thm}

\smallskip

In an effort to justify their conjecture, KTT computed a number of examples of
stretched Kostka polynomials for small~$k$.  In \cite[Conj.~3.3]{KTT04},
they conjectured an explicit description of $P_k$ for $0\le k\le 3$.  By now, all
$k\le 2$ cases are known and follow from the existing results on LR coefficients
(see below). Recently, Prasad \cite[Thm.~7.6]{Pra26} proved the $k=3$ case when $\la$ is a rectangle. We complete the proof of the remaining $k=3$ case of the conjecture for general partitions $\la,\mu$:

\smallskip

\begin{thm}[{\rm formerly KTT's conjecture}{}]\label{t:k3}
Let $\la,\ts\mu$ be partitions such that $K_{\la,\ts\mu}=3$. Then:
\[f_{\la,\ts\mu}(N)\,=\,2N+1 \qquad\text{or}\qquad f_{\la,\ts\mu}(N)\,=\,\frac{(N+1)(N+2)}{2}\..
\]
\end{thm}

\smallskip

We outline the proof of Theorem~\ref{t:finite} later in this section. It hinges
on a geometric result (Main Lemma~\ref{lem:finite}) which is of independent
interest. Theorem~\ref{t:k3} is a mostly trivial consequence of this proof.

\smallskip
\subsection{LR coefficients}

The KTT conjectures are best understood as a stepping stone towards
more general KTT conjectures for \defn{Littlewood--Richardson {\rm (LR)} coefficients} (see below).
The LR coefficients are defined as the structure constants in the ring of Schur functions:
\[s_\la\cdot s_\mu \ =\ \sum_\nu c^\nu_{\la,\ts\mu}\ts s_\nu\..\]
Here $\la,\mu,\nu$ are integer partitions, and $|\la|+|\mu|=|\nu|$.
These coefficients play a central role in algebraic combinatorics,
representation theory, and enumerative geometry, see e.g.\ \cite{Ful97,Mac95,Sta99}.

First, note that Kostka numbers are a special case of LR coefficients. Indeed, for $\ell:=\ell(\mu)$, define:
\[
\tau_i\,:=\,\mu_i +\ldots + \mu_\ell, \qquad \sigma_i\,:=\,\mu_{i+1} + \ldots + \mu_\ell \qquad (1\le i\le\ell).\]
Then \. $K_{\la,\ts\mu}=c^\tau_{\sigma,\ts\la}$, and this identification commutes with stretching; see \cite[Eq.~(1.5)]{KTT04}.

The vanishing of LR coefficients became prominent following Klyachko's work \cite{Kly98},
which connected it to Horn's problem on the eigenvalues of Hermitian matrices.
Knutson and Tao \cite{KT99} completed the argument by establishing the \defn{saturation property:}
\[(\ast)\qquad c^\nu_{\la,\ts\mu}>0 \quad\Longleftrightarrow\quad c^{N\nu}_{N\la,\ts N\mu}>0 \qquad\text{for every integer }N\ge 1.\]
We refer to \cite{Ful00} for background on Horn's problem and this remarkable connection.

The proof brought discrete geometry to the forefront. Knutson and Tao used a \defn{hive polytope} $\rQ^\nu_{\la,\ts\mu}\subseteq\R^D$ whose lattice points are counted by $c^\nu_{\la,\ts\mu}$, and showed that every nonempty hive polytope contains a lattice point \cite{KT99}.
Here $D=\binom{\ell-1}{2}$ is the size of the hive, where $\ell$ bounds the length of the three boundary partitions.
The hive model is linearly equivalent to the polyhedral model arising from the tableau rule, but has simpler defining inequalities and additional symmetries \cite{PV05}. Stretching the boundary partitions dilates the polytope, so
\[f^\nu_{\la,\ts\mu}(N)\,:=\,c^{N\nu}_{N\la,\ts N\mu}\,=\,\bigl|N\cdot\rQ^\nu_{\la,\ts\mu}\cap\Z^D\bigr|.\]
For the specialization above to Kostka numbers, we write $\rQ_{\la,\ts\mu}:=\rQ^\tau_{\sigma,\ts\la}$.

Since hive polytopes have rational vertices, Ehrhart theory implies that $f^\nu_{\la,\ts\mu}(N)$ is a quasipolynomial; see \cite{Bar97,BR15}. Their vertices need not be integral, so polynomiality does not follow from the usual lattice-polytope case of Ehrhart's theorem. Nevertheless, Derksen and Weyman \cite{DW02} proved that these functions are polynomials, resolving the polynomiality part of \cite[Conj.~3.1]{KTT04}. Rassart \cite{Ras04} gave another proof using vector partition functions and chamber complexes; see also \cite{Tha25} for a recent proof.

\smallskip
\subsection{Prior work and related results}

While the hive polytopes need not have integral vertices, the possible
Ehrhart polynomials that arise appear quite constrained. With this in mind,
KTT conjectured that there are only finitely many polynomials  \ts
$f_{\la,\ts\mu}^\nu$ \ts with a fixed value of \ts $c_{\la,\ts\mu}^\nu=k$, see \cite[Conj.~3.3]{KTT04}.
In particular, for $k\le 3$, KTT conjectured the following complete classification:
\begin{align*}
k=0 &\quad\Longrightarrow\quad f_{\la,\ts\mu}^\nu(N)=0 \.,\\
k=1 &\quad\Longrightarrow\quad f_{\la,\ts\mu}^\nu(N)=1\.,\\
k=2 &\quad\Longrightarrow\quad f_{\la,\ts\mu}^\nu(N)=N+1\.,\\[-7pt]
k=3 &\quad\Longrightarrow\quad f_{\la,\ts\mu}^\nu(N)=2N+1 \quad\text{or} \quad f_{\la,\ts\mu}^\nu(N)=\frac{(N+1)(N+2)}{2}\.,
\end{align*}
see \cite[Conj.~3.3]{KTT04}.  These implications have a curious history of
development.

Notably, the $k=0$ case is equivalent to the saturation property,
proved by Knutson and Tao \cite{KT99} as mentioned above.  The $k=1$ case is known as
\defn{Fulton's conjecture}, proved by  Knutson--Tao--Woodward \cite[\S 6.1]{KTW04}.
The $k=2$ case was proved by Ikenmeyer \cite[Thm.~1.1]{Ike16}; see also \cite{She17}
for an algebraic generalization.

Note that while our Theorem~\ref{t:k3} resolves the $k=3$ case for Kostka
coefficients, it remains open for general LR coefficients.  Similarly,
the general finiteness of the number of LR polynomials with given LR coefficient
remains wide open.


A separate conjecture of KTT asserts that the coefficients of stretched LR polynomials in the monomial basis are nonnegative \cite[Conj.~3.1]{KTT04}; the Kostka specialization appears in \cite[Conj.~3.2]{KTT04}. This would imply the saturation property, but remains open in general, even for Kostka numbers, where the saturation property is trivial. Ferudun \cite{Fer26} gives a computer-assisted proof of this conjecture for partitions of length at most five. We refer to \cite{Ale19} for examples illustrating the limitations of positivity statements for related polytopes.

There is also a computational perspective on small coefficients. The saturation theorem gives a polynomial-time algorithm for deciding whether an LR coefficient vanishes, by testing feasibility of its hive polytope using linear programming \cite{DM06,MNS12}. B\"urgisser and Ikenmeyer \cite{BI13} gave a faster combinatorial algorithm, while Ikenmeyer \cite{Ike16} showed that, for every fixed $k\ge 0$, one can decide in polynomial time whether $c^\nu_{\la,\ts\mu}=k$.  The second author and Robichaux gave a probabilistic polynomial-time algorithm for the vanishing
of (more general) Schubert coefficients \cite{PR25}.

Finally, polynomiality is a special feature of the Kostka and LR settings, rather than a general consequence of stretching representation-theoretic multiplicities. Analogous polynomiality statements fail for tensor product multiplicities in other root systems \cite{DM06} and for Newell--Littlewood numbers \cite{GOY21}. For a generalization of stretching to Schubert coefficients, the second and third authors proved eventual quasipolynomiality \cite[Thm.~1.1]{PS26}, and Alexandersson subsequently proved eventual polynomiality \cite[Cor.~7.3]{Ale26}. The qualification ``eventual'' is necessary in that setting \cite[Thm.~1.3]{PS26}.  See also \cite{PR27}
for an overview of the saturation property for stretched Schubert coefficients and other generalizations of LR coefficients.

\smallskip
\subsection{Proof outline}

In our proof of Theorem~\ref{t:finite}, we show that it suffices to consider only partitions with bounded length $\ell(\lambda), \ell(\mu) \le k+2$. To do this, we factor the stretched Kostka polynomial into \emph{primitive Kostka polynomials} coming from \emph{primitive pairs} of partitions $\lambda,\mu$ and use McAllister's degree formula \cite[Thm. 1.1]{McA08}. We will also need the following result on $h^\ast$-positivity of the stretched Kostka polynomials which follows easily from \cite[Props.~3.3 and 3.5]{HMcA11}; see Section~\ref{s:hstar}.

\begin{lemma}[{\rm $h^\ast$-positivity Lemma}{}]\label{l:hstar}
Let $\la,\ts\mu$ be partitions such that $K_{\la,\ts\mu}\ge 1$. Let $d:=\deg(f_{\la,\ts\mu})$, and write the corresponding Ehrhart series as
\begin{equation}\label{eq:Et1}
    \Ehr_{\la,\ts\mu}(t) \,:=\, \sum_{N=\ts0}^\infty f_{\la,\ts\mu}(N)\ts t^N \,=\, \frac{h^\ast_{\la,\ts\mu}(t)}{(1-t)^{d+1}}\..
\end{equation}
Then,
\[h^\ast_{\la,\ts\mu}(t)\ts\in\ts\nn[t].\]
\end{lemma}

\smallskip

We remark that Lemma~\ref{l:hstar} resolves another KTT conjecture \cite[Conj.~4.2]{KTT06}. To complete the proof, we consider Gelfand--Tsetlin polytopes whose Ehrhart polynomials are stretched Kostka polynomials. We then induct on the dimension of the polytope to prove the following discrete geometry lemma which is of independent interest.

\begin{lemma}[{\rm Main Lemma}{}]\label{lem:finite}
Let $k$ be a positive integer, and let \. $A\in \zz^{n\times d}$ \. be an integral $n\times d$ matrix.
Then, up to integer translation, there are at most \ts $\ga=\ga(A,k)$ \ts polytopes of the form \. $Q_\bb := \{\bx \in \R^d: A\bx\le \bb\}$ \.
for \. $\b\in \Z^n$ such that \. $1\le |Q_\bb \cap \zz^d|\le k$.
\end{lemma}

\smallskip

We note that there is ample literature on \emph{lattice point free} \ts
polytopes with given normals, much of it inspired by Herb Scarf's early work,
see e.g.\ \cite{BHS94,JLR07}.  While there is a great deal of literature
on polytopes with $\le k$ points as well, much of it describes the structure
of $Q_\bb$, see e.g.\ \cite{Stu95}.  We were not able to find any finiteness
results on $\ga(A,k)$ available in the literature.

\smallskip
\subsection{Structure of the paper}

We start with some basic definitions and notation in Section~\ref{s:defns}. In Section~\ref{s:discgeo}, we give the proof of the Main Lemma~\ref{lem:finite}. In Section~\ref{s:GT}, we describe the Gelfand--Tsetlin polytope $\GT_{\lambda,\mu}$ corresponding to partitions $\lambda,\mu$ and in Section~\ref{s:hstar}, we show the $h^\ast$-positivity Lemma~\ref{l:hstar}. We prove Theorem~\ref{t:finite} in Section~\ref{s:prooffinite}, and we show the specific case of Theorem~\ref{t:k3} in Section~\ref{s:proofk3}. We conclude with final remarks and open problems in Section~\ref{s:finalrem}.

\medskip
\section{Basic definitions and notation}\label{s:defns}

\subsection{Partitions and dominance}

We write $\N:=\Z_{\ge 0}$. A \defn{partition} $\la=(\la_1,\ldots,\la_\ell)$ is a weakly decreasing sequence of positive integers. Its \defn{size} is $|\la|:=\la_1+\cdots+\la_\ell$, and its \defn{length} is $\ell(\la):=\ell$. We write $\la\vdash n$ when $|\la|=n$. We identify partitions that differ only by trailing zeros, and set $\la_i:=0$ for $i>\ell(\la)$. In particular, padding a partition with zeros does not change its length.

For partitions $\la,\mu\vdash n$, we write $\la\trianglerighteq\mu$ if $\la$ \defn{dominates} $\mu$, that is,
\[\sum_{i\ts=\ts1}^j\la_i\,\ge\,\sum_{i\ts=\ts1}^j\mu_i \qquad\text{for every} \quad j\ge 1.\]
Note that $\la\trianglerighteq\mu$ implies $\ell(\la)\le\ell(\mu)$. A pair $(\la,\mu)$ with $\la\trianglerighteq\mu$ is called \defn{primitive} if
\[\sum_{i=1}^j\la_i\,>\,\sum_{i=1}^j\mu_i \qquad\text{for all} \quad 1\le j<\ell(\mu).\]
Thus primitivity can be thought of as a strict dominance condition.

We use \defn{multiplicity notation} $\la=(a_1^{v_1},\ldots,a_m^{v_m})$ to indicate that the distinct part $a_p$ occurs $v_p$ times, where $a_1>\cdots>a_m$. When a partition is padded to a specified length, this notation includes the multiplicity of zero.

\subsection{Kostka numbers and stretching}

The Young diagram of $\la$ has $\la_i$ boxes in row $i$, with rows numbered from top to bottom. A \defn{semistandard Young tableau} of shape $\la$ and content $\mu$ is a filling of this diagram with positive integers, weakly increasing along rows and strictly increasing down columns, with exactly $\mu_i$ entries equal to $i$. The \defn{Kostka number} $K_{\la,\ts\mu}$ counts these tableaux. We use the standard characterization
\[K_{\la,\ts\mu}>0 \quad\Longleftrightarrow\quad \la\trianglerighteq\mu,\]
for partitions $\la,\mu$ of the same size; see \cite{Ful97,Mac95,Sta99}.

For an integer $N\ge 1$, let $N\la:=(N\la_1,N\la_2,\ldots)$. As mentioned in the introduction, the function $N\mapsto K_{N\la,\ts N\mu}$ agrees with a polynomial in $\mathbb{Q}[N]$, which we denote by
\[f_{\la,\ts\mu}(N)\,:=\,K_{N\la,\ts N\mu} \qquad(N\ge 1).\]
We call $f_{\la,\ts\mu}$ the \defn{stretched Kostka polynomial}. Its value at $N=0$ is understood as polynomial evaluation. If $K_{\la,\ts\mu}=0$, dominance implies that $f_{\la,\ts\mu}$ is identically zero. If $K_{\la,\ts\mu}>0$, the polyhedral interpretation below gives $f_{\la,\ts\mu}(0)=1$.

\subsection{Polytopes and Ehrhart series}

A \defn{rational polyhedron} in $\R^D$ is an intersection of finitely many closed halfspaces defined by inequalities with rational coefficients. A bounded rational polyhedron is a \defn{rational polytope}, and a \defn{lattice polytope} is a polytope whose vertices belong to $\Z^D$. The dimension of a polytope is the dimension of its affine hull, which may be smaller than the ambient dimension $D$.

For a nonempty rational polytope $P\subseteq\R^D$ and an integer $N\ge 0$, let $N\cdot P:=\{N\x:\x\in P\}$. Its \defn{Ehrhart counting function} is:
\[L_P(N)\,:=\,|N\cdot P\cap\Z^D|.\]
Ehrhart's theorem states that $L_P$ is a quasipolynomial of degree $\dim(P)$, and is a polynomial when $P$ is a lattice polytope; see \cite{BR15}. Here, we say a function $f(N)$ is a \defn{quasipolynomial} if there exist an integer $m$ and polynomials $g_i$ for $0\le i \le m-1$ such that $f(N)=g_i(N)$  for all $N=i \! \mod m$. Rational polytopes with nonintegral vertices can also have polynomial Ehrhart counting functions. This is known as period collapse.

Translation by a lattice vector preserves the Ehrhart counting function: if $P'=P+\mathbf z$ for $\mathbf z\in\Z^D$, then $N\cdot P'=N\cdot P+N\mathbf z$, so $L_{P'}(N)=L_P(N)$.

The Gelfand--Tsetlin polytopes described in Section~\ref{s:GT} realize stretched Kostka polynomials as Ehrhart counting functions \cite{McA08}. For $K_{\la,\ts\mu}>0$, set $d:=\deg(f_{\la,\ts\mu})$ and write the corresponding \defn{Ehrhart series} as
\[\Ehr_{\la,\ts\mu}(t)\,:=\,\sum_{N\ts =\ts 0}^{\infty}f_{\la,\ts\mu}(N)\ts t^N\,=\,\frac{h^\ast_{\la,\ts\mu}(t)}{(1-t)^{d+1}}.\]
The numerator $h^\ast_{\la,\ts\mu}(t)\in\Z[t]$ has degree at most $d$. Writing
\[h^\ast_{\la,\ts\mu}(t)\,=\,\sum_{i\ts=\ts0}^d h_i^*t^i,\]
we call $(h_0^*,\ldots,h_d^*)$ the corresponding \defn{$h^\ast$-vector}, allowing trailing zeros. Equivalently,
\begin{equation}\label{eq:binomexp}
    f_{\la,\ts\mu}(N)\,=\,\sum_{i\ts=\ts0}^d h_i^*\binom{N+d-i}{d}.
\end{equation}
In particular, $h_0^*=f_{\la,\ts\mu}(0)=1$.

\medskip
\section{Auxiliary discrete-geometry lemmas}\label{s:discgeo}

\subsection{\texorpdfstring{$A$}{A}-polytopes}

For a matrix $A\in \Z^{n\times d}$ with rows $\a_1,\ldots, \a_n$ and a vector $\b=(b_1,\ldots, b_n)\in \Z^n$, define the polyhedron:
\[P_A(\b) \, := \, \{\x \in \R^d : A\ts\x \, \le \, \b\},\]
where we consider the inequality as a series of componentwise conditions, i.e.,
\[A\ts\x \, \le \, \b \quad \Longleftrightarrow \quad \a_i\x \, \le \, b_i \text{ for } 1\le i \le n.\]
For $1\le i \le n$ and $b_i\in \Z$, define the halfspace:
\[H_i(b_i)\, :=\, \{\x \in \R^d: \a_i \x \le b_i\},\]
so
\[P_A(\b) \, = \, \bigcap_{i\ts=1}^n H_i(b_i).\]

If $\b\in \Z_{\ge 0}^n$ and $P_A(\b)$ is bounded, we call it an $A$\emph{-polytope}. Note that every such polytope contains the origin. Our goal is to prove that, for a fixed matrix $A$, there are only finitely many equivalence classes of $A$-polytopes with few lattice points under the following relation.

We call two $A$-polytopes \emph{equivalent} if one is obtained from the other by translation by a vector in $\Z^d$. Our first tool is a procedure that shrinks an $A$-polytope as much as possible without changing its lattice points.

We say an $A$-polytope $P_A(\b)$ is \emph{reduced} if, for every $i$, the hyperplane $\{\x \in \R^d : \a_i\x=b_i\}$
contains a lattice point of $P_A(\b)$. Given an $A$-polytope $P_A(\b)$, its \emph{reduction} is given by $P_A(\b')$ for $\b':=(b_1',\ldots, b_n')$ where
\[b_i' \, := \, \max\ts\left\{\a_i\x:\x \in P_A(\b)\cap \Z^d\right\}\nts\nts.\]
It is easy to see that $P_A(\b')$ is reduced and satisfies
\[P_A(\b)\cap\Z^d=P_A(\b')\cap\Z^d.\]

\begin{ex}
Figure~\ref{pict:a-polytopes} shows four $A$-polygons in $\R^2$ for
\[A \, = \, \begin{pmatrix}
    -1 & 0\\
    0 & -1\\
    2 & 2\\
    -1 & 2
\end{pmatrix}\nts\nts\nts.\]
From left to right, the polygons are defined by the vectors:
\[\b=(0,\ts0,\ts2,\ts2),\quad \b=(2,\ts0,\ts2,\ts1),\quad \b=(0,\ts0,\ts3,\ts2),\quad\text{and}\quad \b=(0,\ts1,\ts1,\ts0), \]
respectively. In particular, the second polygon is defined by:
\[ -x_1\le 2,\qquad -x_2\le 0,\qquad 2x_1+2x_2\le 2,\qquad -x_1+2x_2\le 1.\]

The third and fourth polygons are equivalent because they are lattice translates of one another. Only the first polygon is reduced. For each of the other polygons, one can decrease an entry of $\b$ without changing its lattice points. In particular, the first polygon is the reduction of the third. The reduction of the second polygon is generated by $\b'=(1,\ts0,\ts2,\ts1)$ and the reduction of the fourth is generated by $(0,\ts1,\ts0,\ts0)$.

\begin{figure}
\centering
\includegraphics[width=0.18\textwidth]{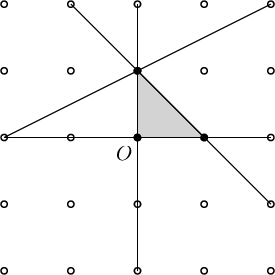}\hskip1.2cm
\includegraphics[width=0.18\textwidth]{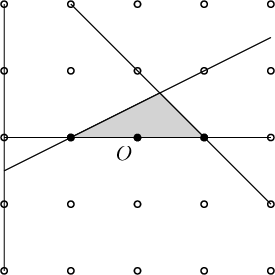}\hskip1.2cm
\includegraphics[width=0.18\textwidth]{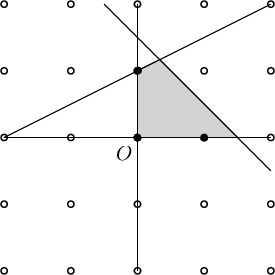}\hskip1.2cm
\includegraphics[width=0.18\textwidth]{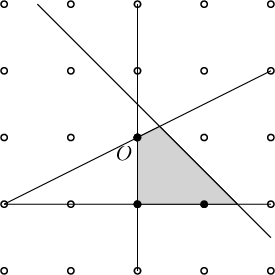}
\caption{Four $A$-polygons.}
\label{pict:a-polytopes}
\end{figure}
\end{ex}

\smallskip
\subsection{Proof of Main Lemma~\ref{lem:finite}}
We first prove that a fixed reduced $A$-polytope can be the reduction of only finitely many $A$-polytopes, up to equivalence.

\begin{lemma}\label{lem:reduced}
Fix a matrix $A\in \Z^{n\times d}$ and let $P_A(\b_0)$ be a reduced $A$-polytope. Then there are only finitely many equivalence classes of $A$-polytopes whose reduction is $P_A(\b_0)$.
\end{lemma}

\begin{proof}
Let $\a_1,\ldots, \a_n$ be the rows of $A$. Suppose, for a contradiction, that there is an infinite collection $\mathcal B\subset \Z_{\ge 0}^n$ of vectors generating pairwise inequivalent $A$-polytopes with common reduction $P_A(\b_0)$. We iteratively construct a set $W\subseteq\{1,\ldots, n\}$ and a tuple of integers $B=(b_i')_{i\ts\in\ts W}$ such that
\[\mathcal B_{W,B}:=\bigl\{\b=(b_1,\ldots, b_n)\in\mathcal B:b_i=b_i'\text{ for every }i\in W\bigr\}\]
is infinite. Initially, $W$ is empty so $\mathcal B_{W,B} = \mathcal B$ is infinite by assumption.

At each stage, let
\[P_{W,B}:=\bigcap_{i\ts\in\ts W}H_i(b_i').\]
If $P_{W,B}$ is bounded, we stop the process. If not, $P_{W,B}$ is an unbounded rational polyhedron. Moreover, $P_{W,B}$ contains the origin so must contain infinitely many lattice points. In particular, there is a lattice point $\mathbf y=(y_1,\ldots,y_d)\in P_{W,B}\setminus P_A(\b_0)$.

For each $\b=(b_1,\ldots, b_n)\in \mathcal B_{W,B}$, we have
\[P_A(\b)\cap \Z^d \, = \, P_A(\b_0)\cap \Z^d.\]
Thus, $\mathbf y\notin P_A(\b)$ so there is at least one index $j\in \{1,\ldots, n\}$ depending on $\b$ such that $\mathbf y\notin H_j(b_j)$, i.e., $\a_j \mathbf y>b_j$. We have $j\notin W$ since $\mathbf y\in P_{W,B}$ and $b_i=b_i'$ for every $i\in W$. Note that $b_j\ge 0$ as $\mathcal{B}_{W,B}\subset \Z_{\ge 0}^n$. Since $\mathcal B_{W,B}$ is infinite, there is thus some index $j\in \{1,\ldots, n\}$ such that the following set is infinite:
\[\left\{\b=(b_1,\ldots,b_n) \in \mathcal{B}_{W,B} : 0\le b_j < \a_j \mathbf y\right\}\nts\nts.\]
But there are only finitely many possibilities for $b_j$ so there is some $b_j'\in \{0,1,\ldots, \a_j\mathbf y -1\}$ such that
\[\mathcal B_{\ts W\ts\cup \ts\{j\}, \ts\ts B\ts\oplus\ts (b_j')}\]
is infinite. Here we use $\oplus$ to indicate concatenation of tuples. We then continue the process.

The process ends either with $W=\{1,\ldots, n\}$ or with $P_{W,B}$ bounded. The first case is impossible because all coordinates of $\b$ are then fixed, so $\mathcal B_{W,B}$ consists of at most one vector. In the second case, for each $j \in \{1,\ldots, n\}\setminus W$, there are only finitely many choices of $b_j'\in \Z_{\ge 0}$ such that the hyperplane $\{\x \in \R^d: \a_j\x=b_j'\}$ meets the bounded polyhedron $P_{W,B}$. If the hyperplane does not meet $P_{W,B}$, then $0\in P_{W,B}$ and $b_j'\ge 0$ imply that all of $P_{W,B}$ lies in the corresponding halfspace so this inequality is redundant. Thus only finitely many pairwise inequivalent polytopes can occur, again a contradiction.
\end{proof}

We now prove Lemma~\ref{lem:finite}.

\begin{proof}
Let $\b\in\Z^n$ satisfy $1\le |P_A(\b)\cap \Z^d|\le k$. Choose a lattice point $\x_0\in P_A(\b)$ and define
\[\b'\,:=\,\b-A\ts \x_0.\]
Since $\x_0\in P_A(\b)$, we have $A\ts \x_0\le \b$ so $\b'\in\Z_{\ge 0}^n$. Moreover,
\[\x\in P_A(\b) \quad\Longleftrightarrow\quad A\ts\x \le \b \quad\Longleftrightarrow\quad  A(\x -\x_0)\le \b' \quad\Longleftrightarrow\quad \x-\x_0\in P_A(\b').\]
Thus,
\[P_A(\b')=P_A(\b)-\x_0,\]
so $P_A(\b)$ is an integer translate of the $A$-polytope $P_A(\b')$. We may therefore assume $\b\in\Z_{\ge 0}^n$.

Discard any zero rows of $A$, whose inequalities are redundant because $\b\ge 0$. Let $\a_1,\ldots, \a_n$ be the remaining rows of $A$. We proceed by induction on the ambient dimension $d$. By Lemma~\ref{lem:reduced}, it is enough to consider reduced $A$-polytopes. The cases $d=0$ and $d=1$ are immediate. Assume the statement holds in every dimension smaller than $d$.

Fix $i\in \{1,\ldots, n\}$. A reduced $A$-polytope contained in the hyperplane $\{\x\in \R^d:\a_i\x=0\}$ can be regarded as a polytope in the rank-$(d-1)$ lattice $\a_i^\perp\cap\Z^d$. The other inequalities have a fixed finite set of integral normals in this lattice, after inequalities parallel to $\a_i$ are discarded. The induction hypothesis therefore gives only finitely many equivalence classes with at most $k$ lattice points. After translating a supporting hyperplane $\{\x \in \R^d: \a_i\x=b_i\}$ by one of its lattice points, the same conclusion applies there. Every lower-dimensional polytope is contained in at least one of its defining hyperplanes, so it remains to consider full-dimensional reduced $A$-polytopes.

For a full-dimensional reduced $A$-polytope $P_A(\b)$, let $\mathbf v(\b)\in\R^n$ record the $(d-1)$-dimensional volumes of the intersections:
\[P_A(\b)\cap\{\x\in \R^d:\a_i\x=b_i\},
\qquad 1\le i\le n,\]
with the convention that the volume is zero when the intersection is lower-dimensional. Since the polytope is reduced, every such intersection contains a lattice point. Translate that point to the origin and identify $\a_i^\perp\cap\Z^d$ with $\Z^{d-1}$. Each intersection is then defined by the fixed restrictions of the rows of $A$, and has between $1$ and $k$ lattice points.

By the induction hypothesis, each coordinate of $\mathbf v(\b)$ belongs to a finite set, so only finitely many such vectors occur. If $\mathbf v(\b_1)=\mathbf v(\b_2)$, the Minkowski uniqueness theorem for convex polytopes \cite[Thm.~18.2]{Gru07} (see also \cite[Thm.~36.2]{Pak10}) shows that $P_A(\b_1)$ and $P_A(\b_2)$ are translates since they have the same facet-normal directions and the corresponding facets have the same volumes.

For a fixed translation class, only finitely many inequivalent translates can be reduced $A$-polytopes. Indeed, if $P_A(\b_2)=P_A(\b_1)+\mathbf t$ and both polytopes are reduced, then
\[b_{r,i}=\max_{\x\in P_A(\b_r)}\a_i\x \qquad (r=1,2).\]
Hence $A\mathbf t=\b_2-\b_1\in\Z^n$. Boundedness implies that the rows of $A$ span $\R^d$, and therefore,
\[L_A:=\{\mathbf t\in\R^d:A\mathbf t\in\Z^n\}\]
is a lattice containing $\Z^d$ with finite index. Thus $L_A/\Z^d$ is finite. Consequently, each possible vector $\mathbf v(\b)$ gives only finitely many equivalence classes, completing the induction.
\end{proof}

\medskip
\section{Gelfand--Tsetlin patterns and \texorpdfstring{$A$}{A}-polytopes}\label{s:GT}

In this section, we describe the Gelfand--Tsetlin polytope $\GT^{\ts\ell}_{\la,\ts\mu}$ and construct a matrix $A_\ell\in \Z^{n\times d}$ and a vector $\b_{\la,\ts\mu}^\ell=(b_1,\ldots, b_n)$ such that:
\[P_{A_\ell}(\b_{\la,\ts\mu}^\ell)=\GT^{\ts\ell}_{\la,\ts\mu}.\]

\subsection{Gelfand--Tsetlin polytopes}

Fix a positive integer $\ell$, and let $\la,\mu$ be partitions with $|\la|=|\mu|$ and $\ell(\la),\ts\ell(\mu)\le\ell$. For convenience, set $\la_i=0$ for $i>\ell(\la)$ and similarly set $\mu_i=0$ for $i>\ell(\mu)$.

Let $X_\ell$ be the space of triangular arrays $(x_{i,\ts j})_{1\ts\le \ts i\ts\le\ts j\ts\le\ts\ell}$ with entries in $\R$. We display such an array as
\[\begin{array}{ccccccccc}
x_{1,\ts\ell} && x_{2,\ts\ell} && x_{3,\ts\ell} && \cdots && x_{\ell,\ts\ell}\\
& x_{1,\ts\ell-1} && x_{2,\ts\ell-1} && \cdots && x_{\ell-1,\ts\ell-1}\\
&& \ddots && \cdots && \iddots\\
&&& x_{1,\ts2} && x_{2,\ts2}\\
&&&& x_{1,\ts1}
\end{array}\]

Set $d:=\binom{\ell+1}{2}$, so $X_\ell=\R^d$. Define the polytope $\GT^{\ts\ell}_{\la,\ts\mu}\subseteq X_\ell$ to consist of the arrays $(x_{i,j})\in X_\ell$ satisfying:
\begin{align}
x_{i+1,\ts j+1}\le x_{i,\ts j}\le x_{i,\ts j+1}
&\qquad\text{for }1\le i\le j\le\ell-1,\label{eq:GT1}\\
x_{i,\ts\ell}=\la_i
&\qquad\text{for }1\le i\le\ell,\label{eq:GT2}\\
\sum_{i\ts=1}^j x_{i,\ts j}=\sum_{i\ts=1}^j \mu_i
&\qquad\text{for }1\le j\le\ell.\label{eq:GT3}
\end{align}

The equalities in \eqref{eq:GT2} and \eqref{eq:GT3} are equivalent to the pairs of inequalities:
\begin{align}
x_{i,\ts\ell}\le\la_i
\quad\text{and}\quad
-x_{i,\ts\ell}\le-\la_i
&\qquad\text{for }1\le i\le\ell,\label{eq:GT4}\\
\sum_{i\ts=1}^j x_{i,\ts j}\le\sum_{i\ts=1}^j\mu_i
\quad\text{and}\quad
-\sum_{i\ts=1}^j x_{i,\ts j}\le-\sum_{i\ts=1}^j\mu_i
&\qquad\text{for }1\le j\le\ell.\label{eq:GT5}
\end{align}

As described, for example, in \cite[\S 2]{McA08}, we have
\[f_{\la,\ts\mu}(N) \,=\, K_{N\la,\ts N\mu} \,=\, \bigl|N\cdot \GT^{\ts\ell}_{\la,\ts\mu}\cap\Z^d\bigr|\qquad (N\ge 1).\]

\subsection{Stretched Kostka polynomials with bounded length}

We are now ready to prove the following Lemma.

\begin{lemma}\label{l:lengthfinite}
Fix positive integers $k$ and $\ell$. The following family of polynomials is finite:
\[P_{k,\ts\ell} \,:=\, \bigl\{f_{\la,\ts\mu}(N):K_{\la,\ts\mu}\le k \. \text{ and } \. \ell(\la),\ts\ell(\mu)\le\ell\bigr\}.\]
\end{lemma}

\begin{proof}
For $1\le i\le j\le\ell$, let $\mathbf e_{i,\ts j}\in X_\ell$ be the elementary basis row vector with $1$ in position $(i,\ts j)$ and $0$ elsewhere. Define
\begin{align*}
\mathcal B_\ell &\,:=\, \bigl(\mathbf e_{i+1,\ts j+1}-\mathbf e_{i,\ts j}\bigr)_{1\ts\le\ts i\ts\le\ts j\le\ts\ell-1} \oplus \bigl(\mathbf e_{i,\ts j}-\mathbf e_{i,\ts j+1}\bigr)_{1\ts\le\ts i\ts\le\ts j\ts\le\ts\ell-1}\ts,\\
\mathcal C_\ell &\,:=\, \bigl(\mathbf e_{i,\ts\ell}\bigr)_{1\ts\le\ts i\ts\le\ts\ell} \oplus \bigl(-\mathbf e_{i,\ts\ell}\bigr)_{1\ts\le\ts i\ts\le\ts\ell}\ts,\\
\mathcal D_\ell &\,:=\, \left(\sum_{i\ts=1}^j\mathbf e_{i,\ts j}\right)_{1\ts\le\ts j\ts\le\ts\ell} \oplus \left(-\sum_{i\ts=1}^j\mathbf e_{i,\ts j}\right)_{1\ts\le\ts j\ts\le\ts\ell},\\
\mathcal A_\ell &\,:=\, \mathcal B_\ell\oplus\mathcal C_\ell\oplus\mathcal D_\ell.
\end{align*}
Here, we use $\oplus$ to indicate concatenation of tuples. Let $n:=2\binom{\ell}{2}+2\ell+2\ell=\ell(\ell+3)$ and $d:=\binom{\ell +1}{2}$. Let $A_\ell\in \Z^{n\times d}$ be the matrix whose rows are given by the tuple $\mathcal A_\ell$. The order of the first $2\binom{\ell}{2}$ rows does not matter but the last $4\ell$ rows are ordered as written in the definitions above.

Now let $\la,\ts\mu$ be partitions such that $f_{\la,\ts\mu}\in P_{k,\ell}$. Namely, $K_{\la,\ts\mu}\le k$ and $\ell(\la),\ts\ell(\mu)\le\ell$. If $K_{\la,\ts\mu}=0$, the saturation property~$(\ast)$ implies that $f_{\la,\ts\mu}(N)$ is identically zero, so assume $K_{\la,\ts\mu}\ge 1$. Define $\b_{\la,\ts\mu}^\ell\in \Z^n$ by:
\[\b^{\ell}_{\la,\ts\mu}:=\left(0^{\binom{\ell}{2}}\right)\oplus \left(0^{\binom{\ell}{2}}\right)\oplus(\lambda_i)_{1\ts\le \ts i\ts\le \ell} \ts\oplus (-\lambda_i)_{1\ts\le \ts i\ts\le \ell}\oplus \left(\sum_{i\ts = 1}^j \mu_i\right)_{1\ts\le\ts j \ts\le\ts \ell}\oplus\left(-\sum_{i\ts = 1}^j \mu_i\right)_{1\ts\le\ts j \ts\le\ts \ell}.\]
In this definition, we write $(0^{\binom{\ell}{2}})$ for the zero-tuple of length $\binom{\ell}{2}$ and we interpret $\lambda_i=0$ (resp. $\mu_i=0$) for $i>\ell(\lambda)$ (resp. $i>\ell(\mu)$). Let $\x:=(x_{i,j})_{1\ts\le\ts i\ts\le\ts j\ts \le\ts \ell} \in X_\ell$. By construction, $A_\ell\ts\x\le \b_{\la,\ts\mu}^\ell$ if and only if the $(x_{i,j})$ satisfy conditions \eqref{eq:GT1}, \eqref{eq:GT4}, and \eqref{eq:GT5}. So we have
\[P_{A_\ell}(\b^\ell_{\la,\ts\mu})=\GT^{\ts\ell}_{\la,\ts\mu}.\]
Therefore,
\[f_{\la,\ts\mu}(N) \,=\, \bigl|N\cdot P_{A_\ell}(\b^\ell_{\la,\ts\mu})\cap\Z^d\bigr| \,=:\, f_{A_\ell,\ts\b^\ell_{\la,\ts\mu}}(N).\]
For $\b\in\Z^n$, set $f_{A_\ell,\ts\b}(N):=\bigl|N\cdot P_{A_\ell}(\b)\cap\Z^d\bigr|$ and define
\[P_{A_\ell,\ts k}:= \bigl\{f_{A_\ell,\ts\b}(N): \b\in\Z^n \text{ and } 1\le f_{A_\ell,\ts\b} (1)\le k\bigr\}.\]
Since $1\le K_{\la,\ts\mu}\le k$, we have $f_{A_\ell,\ts\b^\ell_{\la,\ts\mu}}(N)\in P_{A_\ell,\ts k}$.

By Lemma~\ref{lem:finite}, there are only finitely many polytopes $P_{A_\ell}(\b)$ containing between $1$ and $k$ lattice points, up to integer translation. Since $A$-polytopes that are integer translates have the same Ehrhart quasipolynomial, $P_{A_\ell,k}$ is finite. Thus, $P_{k,\ts\ell}\subseteq P_{A_\ell,\ts k}\cup\{0\}$ is also finite.
\end{proof}

\medskip

\section{Proof of \texorpdfstring{$h^\ast$}{hstar}-positivity Lemma~\ref{l:hstar}}\label{s:hstar}

We refer to \cite{HMcA11} for the relevant notation and definitions. In Section 2 of that paper, Howard and McAllister construct a graded weight ring $R:=R(\lambda,\mu)$ whose degree-$N$ homogeneous component has a $\cc$-basis indexed by the SSYTs of shape $N\la$ and content $N\mu$.
Hence,
\begin{equation}\label{eq:dimRN}
    \dim_\cc(R_N)\,=\, K_{N\la,\ts N\mu} \,=\, f_{\la,\ts\mu}(N) \qquad\text{for every }N\ge 0.
\end{equation}

By \cite[Prop.~3.3]{HMcA11}, there is a homogeneous system of parameters $f_1,\ldots, f_s\in R_1$ where $s$ is the Krull dimension of $R.$ Let $S:=\cc[f_1,\ldots, f_s]$ and let $\f$ be the $s$-tuple $(f_1,\ldots, f_s)$. Then, \cite[Prop.~3.5]{HMcA11} implies that $R$ is a free $S$-module so \cite[Thm.~13.37]{MS05} shows that the Hilbert series of $R$ is:
\begin{equation}\label{eq:HRt}
    H(R;t)\,:=\,\sum_{N=\ts0}^\infty\dim_\cc (R_N)\,t^N \,=\, \frac{H(R/\f R;t)}{(1-t)^s}.
\end{equation}

We can describe $H(R/\f R;t)$ as follows. Suppose $R=Sy_1\oplus \cdots \oplus Sy_m$ and let $r:=\max_j\{\deg(y_j)\}$. Then,
\begin{equation}\label{eq:HRfRt}
    H(R/\f R;t)=\sum_{i\ts=\ts0}^rh_it^i, \qquad h_i=\#\{y_j:\deg(y_j)=i\}\in \nn.
\end{equation}

Thus, equations \eqref{eq:dimRN}, \eqref{eq:HRt}, \eqref{eq:HRfRt}, and the definition of $\Ehr_{\la,\ts\mu}(t)$ in \eqref{eq:Et1} imply that
\begin{equation}\label{eq:Et2}
    \Ehr_{\la,\ts\mu}(t)\,=\,\frac{\sum_{i\ts = \ts 0}^r h_it^i}{(1-t)^s},\qquad h_i\in \nn.
\end{equation}
Comparing \eqref{eq:Et1} and \eqref{eq:Et2} completes the proof.  \qed

\medskip
\section{Proof of Theorem~\ref{t:finite}}\label{s:prooffinite}

If $k=0$, then the saturation property~$(\ast)$ ensures that $P_k=\{0\}$ so we may assume $k\ge 1$. Let $\la,\mu$ be partitions with $f_{\la,\ts\mu}(1)=k$. Since $K_{\la,\ts\mu}\ge 1$, we have $|\la|=|\mu|$. By \cite[\S 1]{McA08}, there is a factorization
\[f_{\la,\ts\mu}(N)=\prod_{r\ts=1}^s f_{\la^{(r)},\ts\mu^{(r)}}(N),\]
where, for every $1\le r\le s$, the partitions $\la^{(r)},\ts\mu^{(r)}\vdash n^{(r)}$ satisfy
\begin{equation}\label{eq:strictdom}
\sum_{i\ts=1}^j\la_i^{(r)} > \sum_{i\ts=1}^j\mu_i^{(r)} \qquad\text{for }1\le j<\ell^{(r)}:=\ell(\mu^{(r)}).
\end{equation}
Each positive integer $f_{\la^{(r)},\ts\mu^{(r)}}(1)$ is at most $k$. By Fulton's conjecture, for example \cite[Prop.~20]{Ale15},
\[K_{\la^{(r)},\ts\mu^{(r)}}=1 \quad\Longleftrightarrow\quad f_{\la^{(r)},\ts\mu^{(r)}}(N)=1 \quad\text{for every }N\ge 1.\]
We may therefore discard all factors for which $K_{\la^{(r)},\ts\mu^{(r)}}=1$. For every remaining factor, \eqref{eq:strictdom} implies $\la^{(r)}\trianglerighteq\mu^{(r)}$ in the dominance order, and hence $K_{\la^{(r)},\ts\mu^{(r)}}>1$. We also have $\ell(\la^{(r)})\le\ell(\mu^{(r)})$.

The number of remaining factors is at most $\lfloor\log_2 k\rfloor$ because each contributes a factor of at least $2$ to the product at $N=1$. Consequently, it is enough to prove that the following family of possible nonconstant factors is finite:
\[\widetilde P_k \,:=\, \left\{f_{\la,\ts\mu}(N): 1<K_{\la,\ts\mu}\le k\. \text{ and } \. \sum_{i\ts=1}^j\la_i> \sum_{i\ts=1}^j\mu_i \text{ for }1\le j<\ell(\mu) \right\}\nts\nts.\]

Let $\la,\ts\mu$ index an element of $\widetilde P_k$, set $\ell:=\ell(\mu)$, and set $d:=\deg(f_{\la,\ts\mu})$. Pad $\la$ with zeros to length $\ell$ and write
\[\la=(k_1^{v_1},k_2^{v_2},\ldots, k_m^{v_m})\]
in multiplicity notation. Since $K_{\la,\ts\mu}>1$, we have $|\la|=|\mu|$. By \cite[Thm.~1.1]{McA08},
\begin{equation}\label{eq:dprim}
d \,=\, \binom{\ell-1}{2}-\sum_{p=1}^m\binom{v_p}{2}.
\end{equation}

By Lemma~\ref{l:hstar}, evaluating \eqref{eq:binomexp} at $N=1$ using $h_0^*=1$ gives:
\[k \,\ge\, f_{\la,\ts\mu}(1) \,=\, h_0^*\binom{d+1}{d}+h_1^*\binom{d}{d} \,\ge \, d+1.\]
Moreover, $d>0$ because $K_{\la,\ts\mu}>1$. Formula~\eqref{eq:dprim} then gives $d\ge\ell-3$, and hence
\[3\le\ell\le d+3\le k+2.\]
Therefore,
\[\widetilde P_k\subseteq\widehat P_k \,:=\, \bigl\{f_{\la,\ts\mu}(N):K_{\la,\ts\mu}\le k \text{ and } \ell(\la),\ts\ell(\mu)\le k+2\bigr\}.\]
The family $\widehat P_k=P_{k,k+2}$ is finite by Lemma~\ref{l:lengthfinite}, completing the proof. \qed

\medskip
\section{Proof of Theorem~\ref{t:k3}}\label{s:proofk3}

Set $\ell:=\ell(\mu)$ and pad $\la$ with zeros to length $\ell$. As in Section~\ref{s:prooffinite}, there is a factorization
\[f_{\la,\ts\mu}(N)=\prod_{r\ts=1}^s f_{\la^{(r)},\ts\mu^{(r)}}(N),\]
where the partitions $\la^{(r)},\ts\mu^{(r)}\vdash n^{(r)}$ satisfy \eqref{eq:strictdom}. Since the product of $f_{\la^{(r)},\ts\mu^{(r)}}(1)$ is $K_{\la,\ts\mu}=3$,  exactly one factor has value $3$ at $N=1$, say $f_{\la^{(1)},\ts\mu^{(1)}}(1)=3$. All other factors are identically $1$ by Fulton's conjecture \cite[Prop.~20]{Ale15}. Thus $f_{\la,\ts\mu}(N)=f_{\la^{(1)},\ts\mu^{(1)}}(N)$, and we may assume that
\begin{equation}\label{eq:prim3}
    \sum_{i\ts=1}^j\la_i>\sum_{i\ts=1}^j\mu_i \qquad\text{for }1\le j<\ell.
\end{equation}

As in Section~\ref{s:prooffinite}, we have $3\le\ell\le K_{\la,\ts\mu}+2 = 5$. If $\ell\le 4$, \cite[Thm.~4.8]{KTT04} states that the corresponding hive polytope $\rQ_{\la,\ts\mu}$ has integral vertices. Since $K_{\la,\ts\mu}=3$, the polytope is either a line segment with one interior lattice point or a lattice triangle with no other lattice points. Their Ehrhart polynomials are, respectively,
\[2N+1 \qquad\text{and}\qquad \frac{(N+1)(N+2)}{2}\ts.\]
We may therefore assume $\ell=5$.

Suppose first that $\la_5>0$. The strict dominance condition \eqref{eq:prim3} for $j=4$, together with $|\la|=|\mu|$, gives $\la_5<\mu_5$. Define
\[\la':=(\la_1-\la_5,\ldots,\la_4-\la_5), \qquad \mu':=(\mu_1-\la_5,\ldots,\mu_5-\la_5).\]
Then $f_{\la',\ts\mu'}(N)=f_{\la,\ts\mu}(N)$ for every $N\ge 1$. Indeed, in every semistandard Young tableau of shape $N\la$ and content $N\mu$, the leftmost $N\la_5$ boxes of row $i$ must contain $i$, for $1\le i\le 5$. We may therefore assume $\la_5=0$.

Now \eqref{eq:dprim} and $1\le d\le 2$ force $d=2, m=2$, and either $(v_1,\ts v_2)=(3,\ts2)$ or $(v_1,\ts v_2)=(2,\ts3)$. Thus,
\[\la=(\la_1,\ts\la_1,\ts\la_1,\ts0,\ts0) \qquad\text{or}\qquad \la=(\la_1,\ts\la_1,\ts0,\ts0,\ts0).\]
Let $j=\ell(\la)$, so $j=3$ or $j=2$, respectively. Define:
\[\la':=(\la_1,\ldots,\la_{j-1},\la_j-\mu_5), \qquad \mu':=(\mu_1,\ldots,\mu_4).\]
Then $f_{\la',\ts\mu'}(N)=f_{\la,\ts\mu}(N)$ for every $N\ge 1$, because in every semistandard Young tableau of shape $N\la$ and content $N\mu$, the rightmost $N\mu_5$ boxes of row $j$ must contain $5$. We have reduced to the case $\ell(\mu')\le 4$, so \cite[Thm.~4.8]{KTT04} again implies that the corresponding hive polytope is integral. The preceding argument completes the proof. \qed

\medskip
\section{Final remarks and open problems}\label{s:finalrem}

\subsection{Positivity}
The main open problem remains positivity of stretched Kostka and LR polynomials in the monomial basis.
After all these years, there is a case to be negative on this conjecture,
especially given some recent notable counterexamples in the area, see e.g.\
\cite{Fer26a}.  Unfortunately, there seems to be relatively little effort
to exhaustively check these conjectures for potential counterexamples, cf.\
\cite{Ale19} for a notable effort.

\subsection{LR finiteness}
One can ask if the LR version of Theorem~\ref{t:finite} can be obtained along
similar lines.  Our approach does not currently extend, since no analogue of the
Howard--McAllister Lemma~\ref{l:hstar} is known in the LR case. While
the KTT conjectures in the LR case may seem abstract and difficult, the following
simple-looking question is perhaps more approachable:
$$
\text{Is the following set bounded:} \quad L(k) \, :=  \,
\big\{c^{2\nu}_{2\la,\ts 2\mu} \, : \, c^\nu_{\la,\ts \mu}=k\big\} \quad \text{for all \. $k\ge 3$}\,?
$$
Most embarrassingly, we don't know if this holds even for $k=3$.

\bigskip

\subsection*{Acknowledgements}
We thank Per Alexandersson, Jes\'us De~Loera, Colleen Robichaux
and Sylvester Zhang for interesting discussions and helpful comments.
The first author was partially supported by Simons Foundation.

\medskip

\subsection*{AI Disclaimer}
ChatGPT 6 Pro was used to proofread and fix minor typos, as well as to verify and correct references and bibliographic details. No AI was used at any point in the process of preparing the mathematical statements and proofs of this paper.
This is because it was largely written in the Spring of 2026,
before we learned the AI tools.  The publication was delayed due to
misplaced hopes in disproving the positivity conjecture.
Had we written the paper now, most likely we would have
used a variety of AI agents.


\end{document}